\documentclass{amsart}
\usepackage{amsmath,amssymb} 
\usepackage{latexsym}
\usepackage{graphicx}
\usepackage{enumitem}

 \usepackage[usenames,dvipsnames]{pstricks}
 \usepackage{epsfig}
 \usepackage{pst-grad} % For gradients
 \usepackage{pst-plot} % For axes

 \usepackage[usenames,dvipsnames]{pstricks}
 \usepackage{epsfig}
 \usepackage{pst-grad} % For gradients
 \usepackage{pst-plot}

\usepackage{tikz}
\usetikzlibrary{arrows}
\usetikzlibrary{patterns}
\usetikzlibrary{shapes}
\usepackage{mathtools} %scommentare per far numerare solo le equationi che sono chiamate nel testo
\mathtoolsset{showonlyrefs,showmanualtags} %scommentare per far numerare solo le equationi che sono chiamate nel testo
\usepackage{bm}
\usepackage{hyperref}
\usepackage{hyperref}
\hypersetup{
			colorlinks=true,
			linkcolor=blue,
			anchorcolor=black,
			citecolor=blue
}
\newcommand{\R}{\mathbb{R}}
\newcommand{\C}{\mathbb{C}}
\newcommand{\Z}{\mathbb{Z}}

\newcommand{\E}{\mathbb{E}\,}

\newcommand{\RP}{\mathbb{R}\textrm{P}}
\newcommand{\CP}{\mathbb{C}\mathrm{P}}

\newcommand{\G}{\mathbb{G}}

\newcommand{\txt}[1]{\textrm{\normalfont{#1}}}
\newcommand{\Oh}{\mathcal{O}}

\newcommand{\be}{\begin{equation}}
\newcommand{\ee}{\end{equation}}

\newcommand{\edeg}{\mathrm{edeg}\,}

\newtheorem{thm}{Theorem}
\newtheorem{lemma}[thm]{Lemma}
\newtheorem{cor}[thm]{Corollary}
\newtheorem{prop}[thm]{Proposition}

\newtheorem*{propi}{Proposition}

\newtheorem*{thmi}{Theorem}

\newtheorem{defi}[thm]{Definition}

\theoremstyle{remark} 
\newtheorem{remark}[thm]{Remark}
\newtheorem{example}[thm]{Example}

\numberwithin{equation}{section}
\numberwithin{thm}{section}
\title{On the number of flats tangent to convex hypersurfaces in random position}
\date{\today}

\author{Khazhgali Kozhasov}
\address{SISSA (Trieste)}
\email{kkozhasov@sissa.it}

\author{Antonio Lerario}
\address{SISSA (Trieste)}
\email{lerario@sissa.it}


 %65H10, 14Q20, 65F22, 65F35, 28A75} 

\begin{document}

\begin{abstract}Motivated by questions in real enumerative geometry \cite{Borcea, PSC, Megyesi, Megyesi2, SottileTheobaldSpheres, SottileTheobald2005, Sottilesurvey}
%and related applications in computational geometry \cite{enclosing2, stab, enclosing, Theobaldvisibility} 
we investigate the problem of the number of flats simultaneously tangent to several \emph{convex} hypersurfaces in real projective space from a probabilistic point of view. More precisely, we say that smooth convex hypersurfaces $X_1, \ldots, X_{d_{k,n}}\subset \RP^n$, where $d_{k,n}=(k+1)(n-k)$, are in random position if each one of them is randomly translated by elements $g_1, \ldots, g_{{d_{k,n}}}$ sampled independently from the orthogonal group with the uniform distribution. Denoting by $\tau_k(X_1, \ldots, X_{d_{k,n}})$ the average number of $k$-dimensional projective subspaces ($k$-\emph{flats}) which are simultaneously tangent to all the hypersurfaces we prove that
 \be \tau_k(X_1, \ldots, X_{d_{k,n}})={\delta}_{k,n} \cdot \prod_{i=1}^{d_{k,n}}\frac{|\Omega_k(X_i)|}{|\textrm{Sch}(k,n)|},\ee
 where ${\delta}_{k,n}$ is the \emph{expected degree} from \cite{PSC} (the average number of $k$-flats incident to $d_{k,n}$ many random $(n-k-1)$-flats), $|\textrm{Sch}(k,n)|$ is the volume of the Special Schubert variety of $k$-flats meeting a fixed $(n-k-1)$-flat (computed in \cite{PSC}) and $|\Omega_k(X)|$ is the volume of the manifold $\Omega_k(X)\subset \G(k,n)$ of all $k$-flats tangent to $X$. We give a formula for the evaluation of $|\Omega_k(X)|$ in terms of some curvature integral of the embedding $X\hookrightarrow \RP^n$ and we relate it with the classical notion of \emph{intrinsic volumes} of a convex set:
 \be\frac{|\Omega_{k}(\partial C)|}{|\textrm{Sch}(k, n)|}=4|V_{n-k-1}(C)|,\quad k=0, \ldots, n-1.\ee
 As a consequence we prove the universal upper bound:
  \be \tau_k(X_1, \ldots, X_{d_{k,n}})\leq {\delta}_{k, n}\cdot 4^{d_{k,n}}.\ee
  Since the right hand side of this upper bound \emph{does not} depend on the specific choice of the convex hypersurfaces, this is especially interesting because already in the case $k=1, n=3$ for every $m>0$ we can provide examples of smooth convex hypersurfaces $X_1, \ldots, X_4$ such that the intersection $\Omega_1(X_1)\cap\cdots\cap\Omega_1(X_4)\subset \G(1,3)$ is transverse and consists of at least $m$ lines.
  
  Finally, we present analogous results for semialgebraic hypersurfaces (not necessarily convex) satisfying some nondegeneracy assumptions.
\end{abstract}

\maketitle

%\tableofcontents

\section{Introduction}
\subsection{Flats simultaneously tangent to several hypersurfaces}Given $d_{k,n} = (k+1)(n-k)$ projective hypersurfaces $X_1, \ldots, X_{{d_{k,n}}}\subset \RP^n$ a classical problem in enumerative geometry is to determine how many $k$-dimensional projective subspaces of $\RP^n$ (called $k$-\emph{flats}) are simultaneously tangent to $X_1, \ldots, X_{{d_{k,n}}}$. 

Geometrically we can formulate this problem as follows. Let $\G(k,n)$ denote the Grassmannian of $k$-dimensional projective subspaces of $\RP^n$ (note that $d_{k,n} = \dim\G(k,n)$). If $X\subset \RP^n$ is a smooth hypersurface, we denote by $\Omega_k(X)\subset \G(k,n)$ the variety of $k$-tangents to $X$, i.e. the set of $k$-flats that are tangent to $X$ at some point. The number of $k$-flats simultaneously tangent to hypersurfaces $X_1, \ldots, X_{{d_{k,n}}}\subset \RP^n$ equals
\be \#\,  \Omega_k(X_1)\cap\cdots \cap \Omega_{k}(X_{{d_{k,n}}}).\ee
Of course this number depends on the mutual position of the hypersurfaces $X_1, \ldots, X_{{d_{k,n}}}$ in the projective space $\RP^n$. 

In \cite{SottileTheobaldSpheres} F.Sottile and T.Theobald proved that there are at most $3\cdot 2^{n-1}$ real lines tangent to $2n-2$ general spheres in $\R^n$  and they found a configuration of spheres with $3\cdot 2^{n-1}$ common tangent lines. They also studied \cite{SottileTheobald2005} the problem of $k$-flats tangent to $d_{k,n}$ many general quadrics in $\RP^n$ and proved that the ``complex bound" $2^{d_{k,n}}\cdot\deg (\G_\C(k,n))$ can be attained by real quadrics. See also \cite{Borcea, Megyesi,Megyesi2,Sottilesurvey} for other interesting results on real enumerative geometry of tangents.

An exciting point of view comes by adopting a random approach: one asks for the \emph{expected value} for the number of tangents to hypersurfaces in \emph{random position}. We say that the hypersurfaces $X_1, \ldots, X_{{d_{k,n}}}\subset \RP^n$ are in \emph{random position} if each one of them is randomly translated by elements $g_1, \ldots, g_{{d_{k,n}}}$ sampled independently from the orthogonal group $O(n+1)$ endowed with the uniform distribution. The average number $\tau_k(X_1,\dots,X_{d_{k,n}})$ of $k$-flats tangent to $X_1, \ldots, X_{d_{k,n}}\subset \RP^n$ in random position is then given by
\be \tau_k(X_1, \ldots, X_{d_{k,n}}):=\E_{g_1,\dots,g_{d_{k,n}}\in O(n+1)}  \#\,  \Omega_k(g_1X_1)\cap\cdots \cap \Omega_{k}(g_{d_{k,n}}X_{{d_{k,n}}}).\ee
The computation and study of properties of this number is precisely the goal of this paper.

 A special feature of the current paper is that we concentrate on the case when the hypersurfaces are boundaries of convex sets. The results we present, however, hold in higher generality as we discuss in Section \ref{sec:nonconvex}.

\begin{defi}[Convex hypersurface] A subset $C$ of $\mathbb{R}\emph{\textrm{P}}^n$ is called (strictly) convex if $C$ does not intersect some hyperplane $L$ and it is (strictly) convex in the affine chart $\mathbb{R}\emph{\textrm{P}}^n\setminus L\simeq \R^n$. A smooth hypersurface $X \subset \R\textrm{\emph{P}}^n$ is said to be \emph{convex} if it bounds a strictly convex open set of $\R\emph{\textrm{P}}^n$.\end{defi}
\begin{remark}[Spherical versus projective geometry]Our considerations in projective spaces run parallel to what happens on spheres, with just small adaptations. A set $C\subset S^n$ is called (strictly) convex if it is the intersection of a (strictly) convex cone $K\subset \R^{n+1}$ with $S^{n}$. A smooth hypersurface $X\subset S^n$ is said to be \emph{convex} if it bounds a strictly convex open set of $S^n$. For the purposes of enumerative geometry, the notion of flats should be replaced with the one of plane sections of $S^n$. Computations involving volumes and the generalized integral geometry formula also require very small modifications (mostly multiplications by a factor of two) and we leave them to the reader. 
\end{remark}

\subsection{Probabilistic enumerative geometry}
Recently, the second author of the current paper together with P. B\"urgisser \cite{PSC}, have studied the similar problem of determining the average number of $k$-flats that simultaneously intersect ${d_{k,n}}$ many $(n-k-1)$-flats in random position in $\RP^n$. They have called this number the \emph{expected degree} of the real Grassmannian $\G(k,n)$, here denoted by ${\delta}_{k,n}$, and have claimed that this is the key quantity governing questions in random enumerative geometry of flats. (The name comes from the fact that the number of solutions of the analogous problem over the complex numbers coincides with the degree of $\G_{\mathbb{C}}(k,n)$ in the Pl\"ucker embedding. Note however that the notion of expected degree is intrinsic and does not require any embedding.)

For reasons that will become more clear later, it is convenient to introduce the \emph{special Schubert variety}\footnote{Note that in the notation of \cite{PSC} we have $\textrm{Sch}(k,n)=\Sigma(k+1, n+1)$ and ${\delta}_{k,n}=\edeg G(k+1, n+1)$.} $\textrm{Sch}(k,n)\subset \G(k,n)$ consisting of $k$-flats in $\RP^n$ intersecting a fixed $(n-k-1)$-flat. The volume\footnote{Here and below we endow the Grassmannian with the Riemannian metric induced by the spherical Pl\"ucker embedding. The smooth locus of a stratified subset of the Grassmannian inherits a Riemannian metric as well, and the volume is computed with respect to this metric, see Section \ref{sec:volumes} for more details.} of the special Schubert variety is computed in \cite[Theorem 4.2]{PSC} and equals
\be |\textrm{Sch}(k,n)|=|\G(k,n)|\cdot \frac{\Gamma\left(\frac{k+2}{2}\right) \Gamma\left(\frac{n-k+1}{2}\right)}{\Gamma\left(\frac{k+1}{2}\right) \Gamma\left(\frac{n-k}{2}\right)},\ee
where $|\G(k,n)|$ denotes the volume of the Grassmannian (see Section \ref{sec:volumes}). The following theorem relates our main problem to the expected degree (see Theorem \ref{thm:main}).

\begin{thmi}[Probabilistic enumerative geometry]The average number of $k$-flats in $\R{\emph{\textrm{P}}}^n$ simultaneously tangent to convex hypersurfaces  $X_1, \ldots, X_{d_{k,n}}$ in random position equals
 \be \tau_k(X_1, \ldots, X_{d_{k,n}})={\delta}_{k,n} \cdot \prod_{i=1}^{d_{k,n}}\frac{|\Omega_k(X_i)|}{|\emph{\textrm{Sch}}(k,n)|},\ee
 where $|\Omega_k(X)|$ denotes the volume of the manifold of $k$-tangents to $X$.
\end{thmi}

The number ${\delta}_{k,n}$ equals (up to a multiple) the volume of a convex body for which the authors of \cite{PSC} coined the name \emph{Segre zonoid}. Except for ${\delta}_{0,n}={\delta}_{n-1, n}=1$, the exact value of this quantity is not known, but it is possible to compute its asymptotic as $n\to \infty$ for fixed $k$. For example, in the case of the Grassmannian of lines in $\RP^n$ one has \cite[Theorem 6.8]{PSC}
\be\label{eq:edeg1} {\delta}_{1,n}=\frac{8}{3\pi^{5/2}}\cdot\frac{1}{\sqrt{n}}\cdot \left(\frac{\pi^2}{4}\right)^{n}\cdot \left(1+\Oh(n^{-1})\right).\ee
The number ${\delta}_{1,3}$ (the average number of lines meeting four random lines in $\RP^3$) can be written as an integral \cite[Proposition 6.7]{PSC}, whose numerical approximation is ${\delta}_{1,3}=1.7262...$. It is an open problem whether this quantity has a closed formula (possibly in terms of special functions).

This reduces our study to the investigation of the geometry of the manifold of tangents, for which we prove the following result (Propositions \ref{propo 1} and \ref{cor:convex} below).
\begin{propi}[The volume of the manifold of $k$-tangents]For a convex hypersurface $X\subset\R{\emph{\textrm{P}}}^n$ we have
\be
\label{eq:convex} \frac{|\Omega_k(X)|}{|\textrm{\emph{Sch}}(k,n)|}=\frac{\Gamma\left(\frac{k+1}{2}\right)\Gamma\left(\frac{n-k}{2}\right)}{\pi^{\frac{n+1}{2}}}\int_{X}\sigma_k(x) dV_X.
\ee where $\sigma_k:X\to \R$ is the $k$-th elementary symmetric polynomial of the principal curvatures of the embedding $X\hookrightarrow \R\textrm{\emph{P}}^n.$ 
\end{propi}
\begin{remark}After this paper was written it was pointed out to us by P. B\"urgisser that this result can be also derived using a limiting argument from \cite{AmBu}, where the tube neighborhood around $\Omega_k(X)$ is described.
\end{remark}
\begin{example}[Spheres in projective space]\label{Example:spheres} Let $S_{r_i}=\{x_1^2+\dots+x_n^2 = (\tan r_i)^2\, x_0^2\} \subset \RP^n$ be a metric sphere in $\RP^n$ of radius $r_i\in (0,\pi/2),\,i=1,\dots,d_{k,n}$ (see Figure \ref{fig:sphere}). Since all principal curvatures of $S_{r_i}$ are constants equal to $\cot r_i$ and since $|S_{r_i}| = \frac{2\sqrt{\pi}^n}{\Gamma(\frac{n}{2})}(\sin r_i)^{n-1}$ Corollary \ref{cor:convex} gives
%we take $X_i=\{x_1^2+x_2^2+x_3^3=(\tan r)^2x_0^2\}$ (see Figure \ref{fig:sphere}). Then, applying the above results, we get
%\begin{align} \tau_0(X_1, X_2, X_3)&= {\delta}_{0,3}\cdot (2 (\sin r)^2)^3& {\delta}_{0,3}=1,\\
%\tau_1(X_1, X_2, X_3, X_4)&={\delta}_{1,3}\cdot  \left(\frac{8}{\pi} \sin r\cos r\right)^4&{\delta}_{1,3}=1.72...,\\
%\tau_2(X_1, X_2, X_3)&={\delta}_{2,3} \cdot (2 (\cos r)^2)^3&{\delta}_{2,3}=1.\end{align}
%More generally, when $X_i=S_{r_i}=\{x_1^2+\dots+x_n^2 = (\tan r_i)^2\, x_0^2\} \subset \RP^n$ is a $(n-1)$-sphere of radius $r_i\in (0,\frac{\pi}{2})$ in $\RP^n$, all its principal curvatures are constants equal to $\cot r_i$ and Corollary \ref{cor:convex} gives
\be \frac{|\Omega_k(S_{r})|}{|\textrm{Sch}(k,n)|}= \frac{2\Gamma\left(\frac{n+1}{2}\right)}{\Gamma\left(\frac{k+2}{2}\right)\Gamma\left(\frac{n-k+1}{2}\right)}\cdot (\cos r_i)^k (\sin r_i)^{n-k-1},\ee
Combining this into Theorem \ref{thm:main} we obtain
\begin{align}\label{eq:spheres}
\tau_k(S_{r_1},\dots,S_{r_{d_{k,n}}}) 
%&= {\delta}_{k,n}\left(\frac{\Gamma(\frac{k+1}{2})\Gamma(\frac{n-k}{2})}{\sqrt{\pi}^{\,n+1}}{n-1 \choose k}\right)^{d_{k,n}} \prod\limits_{i=1}^{d_{k,n}} \frac{2\sqrt{\pi}^n}{\Gamma(\frac{n}{2})}\cos^k R_i \sin^{n-1-k} R_i\\%
={\delta}_{k,n}\cdot  \prod\limits_{i=1}^{d_{k,n}} \left(\frac{2\Gamma\left(\frac{n+1}{2}\right)}{\Gamma\left(\frac{k+2}{2}\right)\Gamma\left(\frac{n-k+1}{2}\right)}\cdot (\cos r_i)^k (\sin r_i)^{n-k-1}\right).
\end{align}
For a fixed $k$ it is natural to find the maximum of the expectation in the case when all the hypersurfaces are spheres. For example, when $k=1$ one can easily see that $\cos r_i (\sin r_i)^{n-2}$ is maximized at $r_i = \arccos \frac{1}{\sqrt{n-1}} = \frac{\pi}{2}-\frac{1}{n^{1/2}}+O(n^{-1/2})$, which is just a bit smaller than $\frac{\pi}{2}$. Therefore,
\begin{align} \max_{r\in (0, \pi/2)}\frac{|\Omega_k(S_{r})|}{|\textrm{Sch}(k,n)|}&=\frac{4}{\sqrt{\pi}}\cdot \frac{\left(\frac{n-2}{n-1}\right)^{\frac{n-2}{2}}}{(n-1)^{\frac{1}{2}}}\cdot \frac{\Gamma\left(\frac{n+1}{2}\right)}{\Gamma\left(\frac{n}{2}\right)}\\&=\left(\frac{8}{e\pi}\right)^{\frac{1}{2}}\left(1+\frac{1}{2n}+\mathcal{O}(n^{-2})\right).
\end{align}
and, together with \eqref{eq:edeg1} and \eqref{eq:spheres}, this gives
\begin{equation}\label{spherebound}
\begin{aligned} \max_{r_1, \ldots, r_{2n-2}\in (0, \frac{\pi}{2})} \tau_1(S_{r_1}, \ldots, S_{r_{2n-2}})&={\delta}_{1, n}\cdot \left(\left(\frac{8}{e\pi}\right)^{\frac{1}{2}}\left(1+\frac{1}{2n}+\mathcal{O}(n^{-2})\right)\right)^{2n-2}\\
&= \frac{e^2}{3\pi^{\frac{3}{2}}}\cdot \frac{1}{\sqrt{n}}\cdot \left(\frac{2\pi}{e}\right)^{n}\cdot\left(1+\mathcal{O}(n^{-1})\right).
\end{aligned}
\end{equation}

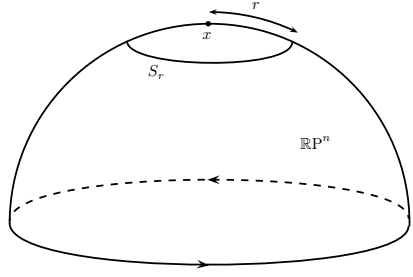
\begin{figure}% Generated with LaTeXDraw 2.0.8
	% Sun Feb 12 18:51:36 CET 2017
	% \usepackage[usenames,dvipsnames]{pstricks}
	% \usepackage{epsfig}
	% \usepackage{pst-grad} % For gradients
	% \usepackage{pst-plot} % For axes
	\scalebox{0.6} % Change this value to rescale the drawing.
	{
		\begin{pspicture}(0,-1.5291992)(8.88,4.7491994)
		\psbezier[linewidth=0.04,linestyle=dashed,dash=0.16cm 0.16cm](0.0,-0.22919922)(0.0,1.0108008)(8.84,0.97080076)(8.84,-0.16919921)
		\psarc[linewidth=0.04](4.43,-0.2791992){4.43}{359.27478}{180.0}
		\psdots[dotsize=0.12](4.4,4.1508007)
		\psbezier[linewidth=0.04](0.0,-0.26919922)(0.0,-1.5091993)(8.84,-1.4691992)(8.84,-0.32919922)
		\psbezier[linewidth=0.04](6.26,3.7393723)(6.26,3.1108007)(2.6,3.1422293)(2.6,3.7708008)
		\psarc[linewidth=0.04,arrowsize=0.05291667cm 2.0,arrowlength=1.4,arrowinset=0.4]{<->}(4.43,-0.019199219){4.43}{63.764862}{90.22646}
		\psline[linewidth=0.072cm,linestyle=dashed,dash=0.16cm 0.16cm,arrowsize=0.05291667cm 2.0,arrowlength=1.4,arrowinset=0.4]{<-}(4.42,-1.2091992)(4.24,-1.1891992)
		\psline[linewidth=0.072cm,linestyle=dashed,dash=0.16cm 0.16cm,arrowsize=0.05291667cm 2.0,arrowlength=1.4,arrowinset=0.4]{<-}(4.38,0.6908008)(4.56,0.6908008)
		\usefont{T1}{ptm}{m}{n}
		\rput(5.461455,4.555801){$r$}
		\usefont{T1}{ptm}{m}{n}
		\rput(4.371455,3.8958008){$x$}
		\usefont{T1}{ptm}{m}{n}
		\rput(6.781455,1.5158008){$\RP^n$}
		\usefont{T1}{ptm}{m}{n}
		\rput(3.251455,3.075801){$S_r$}
		\end{pspicture} 
	}\caption{The equation $x_1^2+\cdots +x_n^2=(\tan r)^2 x_0^2$  defines in $\RP^n$ a metric sphere of radius $r$, i.e. the set of all points at distance $r$ from a fixed point.}\label{fig:sphere}
\end{figure}

We observe that a hypersurface $S_{y,r}$ which is a sphere in some affine chart $U\simeq \R^n$, i.e. $S_{y,r}=\{x\in \R^n: \sum_{i=1}^n (x_i -y_i)^2=r^2\}$, is a convex hypersurface in $\RP^n$, but it is not a sphere with respect to the projective metric unless it's centered at the origin ($y=0$); and, viceversa, a metric sphere in $\RP^n$ needs not be a sphere in an affine chart. In fact, \eqref{spherebound} tells that Sottile and Theobald's upper bound $3\cdot  2^{n-1}$ for the number of lines tangent to $d_{1, n}$ affine spheres in $\R^n$ \emph{does not} apply to the case of spheres in $\RP^n$: since $\frac{2\pi}{e}>2$, when $n$ is large \eqref{spherebound} is larger than $3\cdot 2^{n-1}$; as a consequence there must be a configuration of $d_{1,n}$ projective spheres in $\RP^n$ with (exponentially) more common tangent lines.
\end{example}
\begin{remark}[The semialgebraic case] The theorem above remains true in the case of semialgebraic hypersurfaces $X_1, \ldots, X_{d_{k,n}}\subset \RP^n$ satisfying some mild non-degeneracy conditions (see Section \ref{sec:nonconvex} for more details). Specifically it still holds true that
\be \tau_k(X_1, \ldots, X_{d_{k,n}})={\delta}_{k,n} \cdot \prod_{i=1}^{d_{k,n}}\frac{|\Omega_k(X_i)|}{|\textrm{Sch}(k,n)|},\ee
but the volume of the manifold of $k$-tangents has a more complicated description:
 \begin{align}\frac{|\Omega_k(X)|}{|\textrm{{Sch}}(k,n)|}&=\frac{{n-1\choose k}\Gamma\left(\frac{k+1}{2}\right)\Gamma\left(\frac{n-k}{2}\right)}{\pi^{\frac{n+1}{2}}}\int_{X}\E_{\Lambda\in Gr_k(T_xX)} | B_x(\Lambda)| dV_X,\end{align}
where $|B_x(\Lambda)|$ denotes the absolute value of the determinant of the matrix of the second fundamental form of $X\hookrightarrow \RP^n$ restricted to $\Lambda\in Gr_k(T_{x}X)$ and written in an orthonormal basis of $\Lambda$ (see Section \ref{sec:convex}), and the expectation is taken with respect to the uniform distribution on $Gr_k(T_xX)\simeq Gr(k,n-1).$
\end{remark}

\subsection{Relation with intrinsic volumes} The quantities $|\Omega_k(X)|$ offer an alternative interesting interpretation of the classical notion of \emph{intrinsic volumes}. Recall that if $C$ is a convex set in $\RP^n$, the spherical Steiner's formula \cite[(9)]{GaoSchneider} allows to write the volume of the $\epsilon$-neighborhood $\mathcal{U}_{\,\RP^n}(C,\epsilon)$ of $C$ in $\RP^n$ as
\be |\mathcal{U}_{\,\RP^n}(C,\epsilon)|=|C|+\sum_{k=0}^{n-1}f_{k}(\epsilon)|S^k||S^{n-k-1}|V_{k}(C),\ee
where
\be\label{eq:universal} f_{k}(\epsilon)=\int_{0}^\epsilon (\cos t)^{k}(\sin t)^{n-1-k}dt.\ee The quantities $V_0(C), \ldots, V_{n-1}(C)$ are called intrinsic volumes of $C$.
What is remarkable is that when $C$ is smooth and strictly convex, $|\Omega_k(\partial C)|$ coincides, up to a constant depending on $k$ and $n$ only, with the $(n-k-1)$-th intrinsic volume of $C$ (again this property can be derived by a limiting argument from the results in \cite{AmBu}).
\begin{propi}[The manifold of $k$-tangents and intrinsic volumes]Let $C\subset \R\emph{\textrm{P}}^n$ be a smooth strictly convex set. Then
\be |V_{n-k-1}(C)|=\frac{1}{4}\cdot \frac{|\Omega_{k}(\partial C)|}{|\textrm{\emph{Sch}}(k, n)|},\quad k=0, \ldots, n-1.\ee
\end{propi}
This interpretation offers possible new directions of investigation and allows to prove the following upper bound (see Corollary \ref{cor:upper})
\be\label{eq:uniform} \tau_k(X_1, \ldots, X_{d_{k,n}})\leq {\delta}_{k,n}\cdot 4^{d_{k,n}},\ee
where the right-hand side depends only on $k$ and $n$.
However, already for $n=3$, as observed by T. Theobald there is no upper bound on the number of lines that can be simultaneously tangent to four convex hypersurfaces in $\RP^3$ in general position (see Section \ref{sec:counter} for details).

%\begin{example}\label{example:unbounded}The following example is due to T. Theobald \cite{private}. Let $\ell\subset \RP^3$ be a line and pick \end{example}

\subsection{Related work}Enumerative geometry over the field of complex numbers is classical. Over the Reals it is a much harder subject, due to the nonexistence of generic configurations. From the deterministic point of view we mention, among others, the papers that are closest to our work and that gave a motivation for it: \cite{Borcea, Megyesi, Megyesi2, SottileTheobaldSpheres, SottileTheobald2005, Sottilesurvey}. The probabilistic approach to real enumerative geometry was initiated in \cite{PSC} for what concerns Schubert calculus, and in \cite{BLLP} for the study of the number of real lines on random hypersurfaces.
\subsection*{Acknowledgements}We wish to thank P. B\"urgisser, K. Kohn, F. Sottile and T. Theobald for interesting discussions.
\section{Preliminaries}By $\G(k,n)\simeq Gr(k+1,n+1)$ we denote the Grassmannian of $(k+1)$-planes in $\mathbb{R}^{n+1}$ (or, equivalently, the set of projective $k$-flats in $\RP^n$). Both notations are used throughout the article. The dimension of $\G(k,n)$ is denoted by $d_{k,n} := \dim\G(k,n) = (k+1)(n-k)$.
\subsection{Metrics \& volumes}\label{sec:volumes}
The Grassmannian $Gr(k,n)$ is endowed with an $O(n)$-invariant
riemannian metric through the Pl\"ucker embedding
$$i: Gr(k,n) \hookrightarrow P\left(\bigwedge^k\R^n\right)$$
where $\text{\normalfont{P}}(\bigwedge^k\R^n)$, the projectivization of the vector space $\bigwedge^k\R^n$, is endowed with the standard metric.
Using this we locally identify $Gr(k,n)$ with the set of unit simple
$k$-vectors ${v_1\wedge\cdots\wedge v_k}$, where $v_1,\dots,v_k$ are
orthonormal in $\R^n$ (see \cite{Kozlov} for more details).

A canonical left-invariant metric on the orthogonal group $O(n)$ is defined as
$$\langle A,B \rangle : = \frac{1}{2}\text{\normalfont{tr}}(A^tB),\ A,
B \in T_{\textbf{1}}O(n)$$
Denoting by $|X|$ the total volume of a Riemannian manifold $X$ (whenever it is finite) one can prove the following formulas
\begin{align}
|Gr(k,n)| = \frac{|O(n)|}{|O(k)||O(n-k)|},\quad
\frac{|O(n+1)|}{|O(n)|}=|S^n|,\quad
|O(1)|=2,\quad
|S^n| = \frac{2\pi^{\frac{n+1}{2}}}{\Gamma(\frac{n+1}{2})}.
\end{align}
For an $m$-dimensional semialgebraic subset $X$ of $Gr(k,n)$ by $|X|$ we denote the $m$-dimensional volume of the set of smooth points $X_{sm}$ of $X$.  
\subsection{Probabilistic setup}
Given a Riemannian manifold $Y$ and a smooth function $f: Y\rightarrow \R$ we denote by $\int_Y f(y)\,dV_Y$ 
the integration of $f$ with respect to the Riemannian volume density of $Y$. We recall that there is a unique $O(n)$-invariant probability distribution on $O(n), Gr(k,n)$ and $S^n$ called \emph{uniform} (see \cite{PSC, Kozlov} for more details). For a measurable subset $A \subset X \in \{O(n), Gr(k,n), S^n\}$ it is defined as 
$$\mathbb{P}(A):= \frac{1}{|X|}\int\limits_X \textbf{1}_A\,dV_X.$$
In the sequel all probabilistic concepts are referred to the above listed spaces endowed with the uniform distribution.
\begin{remark}
For a measurable $A \subset Gr(k,n)$ the set $\hat{A} = \{g\in O(n): g^{-1}\mathbb{R}^k \in A\}$ is measurable in $O(n)$ and 
\begin{align*}
\mathbb{P}(A) = \frac{1}{|Gr(k,n)|}\int\limits_{Gr(k,n)} \textbf{1}_{A}\, dV_{Gr(k,n)} = \frac{1}{|O(n)|} \int\limits_{O(n)} \textbf{1}_{\hat{A}}\,dV_{O(n)} = \mathbb{P}(\hat{A})
\end{align*}
We will implicitly use this identification when needed.

\end{remark}
\subsection{Integral geometry of coisotropic hypersurfaces of Grassmannian}
A smooth (respectively semialgebraic) hypersurface $\mathcal{H}$ of $\G(k,n)$ is said to be \emph{coisotropic} if for any (respectively for any smooth point of codimension one) $\Lambda \in \mathcal{H}$ the normal space $N_{\Lambda}\mathcal{H}\subset T_{\Lambda}\G(k,n) \simeq \txt{Hom}(\Lambda,\Lambda^{\perp})$ is spanned by a rank one operator. 

\noindent
For $k,m\geq 1$ let $u_j\in S(\R^k), v_j \in S(\R^m),\, j=1,\dots,km$ be unit independent random vectors. Then the \emph{average scaling factor} $\alpha(k,m)$ is defined as
\begin{align} 
 \alpha(k,m):=\E\Vert (u_1\otimes v_1)\wedge\cdots\wedge(u_{km}\otimes v_{km})\Vert
 \end{align}
where $\Vert\cdot\Vert$ is induced from the standard scalar product on $\R^k\otimes\R^m$: $(u_1\otimes v_1,u_2\otimes v_2):=(u_1,u_2)(v_1,v_2)$.
We will use the generalized Poincar\'e formula for coisotropic hypersurfaces of $\G(k,n)$ proved in \cite[Thm. 3.19]{PSC}:

\begin{thm}\label{Poincare theorem}
Let $\mathcal{H}_1,\dots,\mathcal{H}_{d_{k,n}}$ be coisotropic hypersurfaces of $\G(k,n)$. Then
\begin{align}
  \E \#(g_1\mathcal{H}_1\cap\cdots\cap g_{d_{k,n}}\mathcal{H}_{d_{k,n}})=\alpha(k+1,n-k)\,|\G(k,n)|\,\prod\limits_{i=1}^{d_{k,n}}\frac{|\mathcal{H}_i|}{|\G(k,n)|}
\end{align}
where $g_1,\dots,g_{d_{k,n}} \in O(n+1)$ are independent randomly chosen orthogonal transformations.
\end{thm}
\begin{remark}
This theorem expresses the average number of points in the intersection of $d_{k,n}$ many hypersurfaces of $\G(k,n)$ in \emph{random position} in terms of the volumes of the hypersurfaces and the average scaling factor $\alpha(k+1, n-k)$, which only depends on the pair $(k,n)$.
\end{remark}
\subsection{Intersection of special real Schubert varieties}
A \emph{special real Schubert variety} $\text{\normalfont{Sch}}(k,n)$ consists of all projective $k$-flats in $\RP^n$ that intersect a fixed projective $(n-k-1)$-flat $\Pi$:
$$\text{\normalfont{Sch}}(k,n)=\{\Lambda \in \G(k,n): \Lambda\cap \Pi \neq \varnothing\}$$
It is a coisotropic algebraic hypersurface of $\G(k,n)$.
In \cite{PSC} P. B\"urgisser and the second author of the current article had introduced a notion of \emph{expected degree} $\delta_{k,n}$ of the Grassmannian $\G(k,n)$. It is defined as the average number of projective $k$-flats in $\RP^n$ simultaneously intersecting $d_{k,n}$ many random projective $(n-k-1)$-flats independently chosen in $\G(n-k-1,n)$. In other words,
\begin{align}
\delta_{k,n} := \E \#(g_1\text{\normalfont{Sch}}(k,n)\cap \cdots\cap g_{d_{k,n}}\text{\normalfont{Sch}}(k,n)). 
\end{align}
Using the formula in \cite[Thm. 4.2]{PSC} for the volume of $\text{\normalfont{Sch}}(k,n)$:
\begin{align}
|\text{\normalfont{Sch}}(k,n)|=|\G(k,n)|\frac{\Gamma(\frac{k+2}{2})}{\Gamma(\frac{k+1}{2})}\frac{\Gamma(\frac{n-k+1}{2})}{\Gamma(\frac{n-k}{2})}
\end{align}
and Theorem \ref{Poincare theorem} one can express
\begin{align}
\delta_{k,n} = \alpha(k+1,n-k)\,|\G(k,n)| \left(\frac{\Gamma(\frac{k+2}{2})}{\Gamma(\frac{k+1}{2})}\frac{\Gamma(\frac{n-k+1}{2})}{\Gamma(\frac{n-k}{2})} \right)^{d_{k,n}}
\end{align}
\begin{remark}
The exact value of $\delta_{k,n}$ (equivalently $\alpha(k+1,n-k)$) remains unknown for $0<k<(n-1)$. See \cite[Sect. 6]{PSC} for various asymptotics of $\delta_{k,n}$. 
\end{remark}
\begin{remark}Note that one can define a notion of ``expected degree'' even over the complex numbers, by sampling complex projective subspaces uniformly from the complex Grassmannian.  Denoting by $c_{k,n}\in H^{2}(\G^{\C}(k,n);\Z)$ the first Chern class of the tautological bundle and by $[\G^\C(k,n)]\in H_{2{d_{k,n}}}(\G^\C(k,n);\Z)$ the fundamental class we have that
\be \textrm{the expected degree over the complex numbers}=\left\langle (c_{k,n})^{{d_{k,n}}}, [\G^\C(k,n)]\right \rangle\ee The resulting number also equals the degree of $\G^{\C}(k,n)$ in the Pl\"ucker embedding.
\end{remark}

\section{The manifold of tangents}\label{sec:convex}
Let $X=\partial C$ be a convex hypersurface of $\RP^n$ (bounding the strictly convex open set $C\subset\RP^n$) and let $p: Gr_k(X) \rightarrow X$ be the Grassmannian bundle of $k$-planes of $X$ (this is a smooth fiber bundle over $X$ whose fiber $p^{-1}(x)$ is the Grassmannian $Gr_k(T_xX)\simeq Gr(k,n-1)$). Define the \emph{$k$th Gauss map} 
\begin{align*}
	\psi:\ Gr_k(X) &\rightarrow \G(k,n)\\
	(x, \Lambda)&\mapsto \text{\normalfont{P}}(\textrm{Span}\{x,\Lambda\})
\end{align*}
here we identify the tangent space $T_x\RP^n$ with the hyperplane $x^{\perp} \subset \R^{n+1}$ and thus $\Lambda$ and $x$ (a line in $\R^{n+1}$) are both subspaces of $\R^{n+1}$.

With this notation we observe that $\psi$ is a smooth embedding and that $\Omega_k(X)$, the set of all $k$-flats tangent to $X$, coincides, by definition, with $\textrm{im} (\psi)$. 

Let's choose a unit normal vector field $\nu$ to $X\subset \RP^n$ pointing inside the convex region $C$. Then the second fundamental form $B$ of $X$ is positive definite everywhere. For $(x,\Lambda) \in Gr_k(X)$ and an orthonormal basis $v_1,\dots,v_k$ of $\Lambda$ let's denote by $B_x(\Lambda)=\det( B(v_i,v_j))$ the determinant of the $k\times k$ matrix $\{B(v_i,v_j)\}$. Note that $B_x(\Lambda)$ does not depend on the choice of $v_1,\dots,v_k$. Using the smooth coarea formula we prove the following proposition.

\begin{prop}\label{propo 1}If $X\subset \mathbb{R}\textrm{\emph{P}}^n$ is a convex hypersurface, then
  \begin{align}\label{eq:volume convex}
\left|\Omega_k(X)\right| =  \frac{\left| Gr(k,n-1)\right|}{{n-1 \choose k}}\int\limits_{X} \sigma_k(x) dV_X
\end{align}
where $\sigma_k:X\to \R$ is the $k$-th elementary symmetric polynomial of the principal curvatures of the embedding $X\hookrightarrow \R\textrm{\emph{P}}^n.$
\end{prop}
\begin{proof}
The $O(n+1)$-invariant metric $g$ on $\G(k,n)$ induces a Riemannian metric $\psi^*g$ on $Gr_k(X)$ through the embedding $\psi$. Note that the restriction of $\psi^*g$ to the fibers $Gr_k(T_xX)$ is $O(T_xX)\simeq O(n-1)$-invariant. We apply the smooth coarea formula to $p: (Gr_k(X),\psi^*g) \rightarrow (X,g_X)$, where $g_X$ denotes the induced metric on $X\hookrightarrow\RP^n$. We obtain:
\begin{align}\label{coarea formula}
 |\Omega_k(X)|=\int\limits_{Gr_k(X)} dV_{Gr_k(X)} = \int\limits_{X}\int\limits_{Gr_k(T_xX)} \left(\emph{NJ}_{(x,\Lambda)}p\right)^{-1}dV_{Gr_k(T_xX)}\, dV_X.
\end{align}

Let's show first that the normal Jacobian $\emph{NJ}_{(x,\Lambda)}p$ equals $\left| B_x(\Lambda)\right|^{-1}=\left| \det(B(v_i,v_j))\right|^{-1}$.

Given a point $x \in X$, a unit normal $\nu\in T_x\RP^n$ to $T_xX$ and an orthonormal basis $v_1,\dots,v_k\in T_xX$ of $\Lambda \in Gr_k(T_xX)$ let's complete them to an orthonormal basis $x,\nu,v_1,\dots,v_k,v_{k+1},\dots,v_{n-1}$ of $\R^{n+1}$. Using these vectors we describe the tangent space to $Gr_k(X)$ at $(x,\Lambda)$.

For $i=1,\dots,n-1$ and $j=1,\dots,k$ let $x_i=x_i(t)$ be a small curve through $x$ in the direction $v_i$ and let $v_j^i=v_j^i(t)$ be the parallel transport of $v_j$ along $x_i$, i.e. the vector field solving $\nabla_{\dot{x}_i}^X v_j^i = 0, v_j^i(0)=v_j$. Note that for any time $t$ the vectors $v_1^i(t),\dots,v_k^i(t) \in T_{x_i(t)}X$ remain pairwise orthonormal. Consider now curves in $Gr_k(X)$ and their tangents  produced by these vectors:
\begin{align}
\widetilde{\gamma}_i(t) &= (x_i(t),v_1^i(t)\wedge\cdots\wedge v_k^i(t))\\
\widetilde{\Gamma}_i:&=\dot{\widetilde{\gamma}}_i(0) = (v_i,\sum\limits_{j=1}^k v_1\wedge\cdots\wedge\dot{v}_j^i(0)\wedge\cdots\wedge v_k)
\end{align}
Observe that
$$\dot{v}^i_j(0) = \nabla^{\R^{n+1}}_{v_i} v_j^i = \underbrace{\nabla^X_{v_i} v_j^i}_{=0} + a_{ij}\, x + b_{ij}\,\nu = a_{ij}\, x + b_{ij}\,\nu$$ 
Since the standard scalar product on $\R^{n+1}$ (here denoted by a dot) induces the  metric on $T_x\RP^n = T_xS^n = x^{\perp}$ and since the second fundamental form of the unit sphere $S^n\subset \R^{n+1}$ coincides with the metric tensor we have  
\begin{align}
a_{ij} &= (\nabla^{\R^{n+1}}_{v_i} v_j^i) \cdot x  = \delta_{ij}\\
b_{ij} &= (\nabla^{\R^{n+1}}_{v_i} v_j^i)\cdot \nu = \label{b_{ij}} (\nabla^{\RP^n}_{v_i}v_j^i+\delta_{ij}\, x)\cdot \nu = (\nabla^{\RP^n}_{v_i}v_j^i)\cdot \nu = B(v_i,v_j)
\end{align}
The tangent space to the fiber $T_{(x,\Lambda)}Gr_k(T_xX) = \ker(p_*)$ is spanned by the following $k(n-1-k)$ vectors:
\begin{align}
\widetilde{\theta}_{ij}(t) &= (x,v_1\wedge\cdots\wedge\underset{i}{(v_i\,\cos t+ v_j\,\sin t)}\wedge\cdots\wedge v_k),\ i=1,\dots,k\\ 
\widetilde{\Theta}_{ij} :&= \dot{\widetilde{\theta}}_{ij}(0) = (0,v_1\wedge\cdots\wedge\underset{i}{v_j}\wedge\cdots\wedge v_k),\ j=k+1,\dots,n-1 
\end{align}
%Indeed, $\widetilde{\theta}_{ij}$'s rotate in $T_xX$ and $\widetilde{\Theta}_{ij}$'s are obviously independent.
We work with the images $\Gamma_i,\Theta_{ij}\in T_{\text{\normalfont{Span}}\{x,\Lambda\}}\G(k,n)$ of $\widetilde{\Gamma}_i$ and $\widetilde{\Theta}_{ij}$ under $\psi_*$. It is easy to see that
\begin{align}
\Gamma_i &= \psi_*\widetilde{\Gamma}_i = v_i \wedge v_1\wedge\cdots\wedge v_k + \sum\limits_{j=1}^kb_{ij}\, x \wedge v_1 \wedge \cdots \wedge\underset{j}{\nu} \wedge\cdots \wedge v_k,\quad  1\leq i \leq n-1\\
\Theta_{ij} &= \psi_*\widetilde{\Theta}_{ij} = x\wedge v_1\wedge\cdots\wedge\underset{i}{v_j}\wedge\cdots\wedge v_k,\quad 1\leq i\leq k,\ k+1\leq j\leq n-1
\end{align}
and $\Gamma_i$'s are orthogonal to $\Theta_{ij}$'s, but $\Gamma_i$'s are not in general orthonormal vectors. Therefore, since $p_* \widetilde{\Gamma}_i = v_i$ and the $v_i$'s form an orthonormal basis for $T_xX$ in order to compute the normal Jacobian $NJ_{(x,\Lambda)}p$ we need to find a change of basis matrix from $\{\Gamma_i\}_{1\leq i\leq n-1}$ to some orthonormal basis of $\text{\normalfont{Span}}\{\Gamma_i\}_{1\leq i\leq n-1} = \ker(p_*\circ\psi_*^{-1})^{\perp}$. For this purpose let's note that for the orthonormal vectors

\begin{align}
  S_j &= x\wedge v_1 \wedge \cdots\wedge \underset{j}{\nu} \wedge \cdots\wedge v_k, \quad\quad 1\leq j \leq k\\
  P_i &= v_i \wedge v_1 \wedge \cdots \wedge v_k,\qquad\quad\quad k+1\leq i \leq n-1
\end{align}
we have
\begin{align}
  \begin{pmatrix}
  \Gamma_1 \\  
  \vdots\\
  \Gamma_k\\
  \Gamma_{k+1}\\
  \vdots\\
  \Gamma_{n-1}
  \end{pmatrix}
=
  \begin{pmatrix}
    b & 0\\
    * & 1
  \end{pmatrix}
        \begin{pmatrix}
          S\\R
        \end{pmatrix}
=
  \begin{pmatrix}
  b_{11} & \dots & b_{1k} & 0 & 0 &\dots & 0\\
  \vdots & \vdots &\vdots &\vdots &\vdots &\vdots& \vdots\\
  b_{k1} & \dots & b_{kk} & 0 & 0 & \dots&0\\
  b_{k+1,1} & \dots & b_{k+1,k} & 1 & 0 & \dots & 0\\
  \vdots & \vdots &\vdots &\vdots &\vdots & \vdots&\vdots\\
  b_{n-1,1} & \dots & b_{n-1,k} & 0 & 0 &0 & 1
  \end{pmatrix}
                                             \begin{pmatrix}
                                               S_1\\
                                               \vdots \\
                                               S_k\\
                                               P_{k+1}\\
                                               \vdots\\
                                               P_{n-1}
                                             \end{pmatrix}                                           
\end{align}
where $b=\{b_{ij}\}_{1\leq i,j\leq k} = \{B(v_i,v_j)\}_{1\leq i,j\leq k}$ by \eqref{b_{ij}}. Note that 
\begin{align}\label{criteria}
  \psi_*\ \text{is injective iff}\ b\ \text{is invertible iff}\ B|_{\Lambda}\ \text{is non-degenerate.}
\end{align}
Then since $B$ is positive definite everywhere $b$ is invertible and
\begin{align}\label{new basis}
                                             \begin{pmatrix}
                                               S_1\\
                                               \vdots \\
                                               S_k\\
                                               P_{k+1}\\
                                               \vdots\\
                                               P_{n-1}
                                             \end{pmatrix}                                           
=
  \begin{pmatrix}
    b^{-1} & 0\\
    * & 1
  \end{pmatrix}
  \begin{pmatrix}
  \Gamma_1\\  
  \vdots\\
  \Gamma_k\\
  \Gamma_{k+1}\\
  \vdots\\
  \Gamma_{n-1}
  \end{pmatrix}
\end{align}
Applying $p_*\circ\psi_*^{-1}$ to  the $S_j,P_i$'s we obtain that $$NJ_{(x,\Lambda)}p = |\det(b^{-1})| = |B_x(\Lambda)|^{-1}$$
and thus
\begin{align}\label{eq:step1}
|\Omega_{k}(X)| = \int\limits_X\int\limits_{Gr_k(T_xX)} |B_x(\Lambda)|\,dV_{Gr_k(T_xX)}dV_X.
\end{align}
Since the fibers $Gr_k(T_xX)$ are endowed with $O(n-1)\simeq O(T_xX)\simeq O(\{x,\nu_x\}^{\perp})$-invariant metric we may rewrite the inner integral as
\begin{align}\label{inner integral}
\int\limits_{Gr_k(T_xX)} |B_x(\Lambda)|\,dV_{Gr_k(T_xX)} = |Gr(k,n-1)|\,\E_{\Lambda \in Gr(k,n-1)}|B_x(\Lambda)|
\end{align}
Since the restriction $B|_{\Lambda}$ of a positive definite form $B$ is also positive definite, we have $B_x(\Lambda)>0$ and hence
\be \E_{\Lambda \in Gr(k,n-1)}|B_x(\Lambda)|=\E_{\Lambda \in Gr(k,n-1)}B_x(\Lambda).\ee
We prove that 
\begin{align}\label{expectation of the second fundamental form}\E_{\Lambda\in Gr(k,n-1)}B_x(\Lambda) = {n-1 \choose k}^{-1} s_k(d_1(x),\dots,d_{n-1}(x))
\end{align}
where the $d_i(x)$'s are the principal curvatures of $X \subset \RP^n$ at the point $x$ and $s_k$ is the $k$-th elementary symmetric polynomial.
Now let's choose an o.n.b. $e=\{{\delta}_1,\dots,{\delta}_{n-1}\}$ of $T_xX$ in which the second fundamental form $B$ is diagonal $D=\text{\normalfont{diag}}\{d_1,\dots,d_{n-1}\}$. 
For vectors $v_i$ we denote by the same letters their coordinate representation in the basis $e$. Let $V$ and $E$ be $(n-1) \times k$ matrices with columns $\{v_i\}_{1\leq i\leq k}$ and $\{{\delta}_i\}_{1\leq i\leq k}$ respectively:
\begin{align}
V = 
\begin{pmatrix}
| &  & | \\
v_1 & \dots & v_k \\
| &  & | 
\end{pmatrix} &\qquad
E = 
\begin{pmatrix}
1 & 0 & \dots & 0\\
\vdots & \vdots & \vdots & \vdots \\
0 & \dots &  0 & 1\\
0 & \dots & \dots & 0\\
\vdots & \vdots & \vdots & \vdots \\
0 & \dots & \dots & 0\\
\end{pmatrix}
\end{align}
There exists an orthogonal matrix $g \in O(n-1)$ s.t. $V = g\cdot E$ and then $b=\{B(v_i,v_j)\}_{1\leq i,j\leq k}$ can be written as
$b = V^t D V = E^t  g^t D g E$. In this view $B_x(\Lambda) = \det(b) = \det(E^t  g^t D g E)$ is just the leading principal minor of $g^tDg$ of order $k$. Note that $B_x(\Lambda)$ does not depend on the choice of $g$, namely it's invariant under the action of $\text{\normalfont{Stab}}_{\,\text{\normalfont{Span}}\{{\delta}_1,\dots,{\delta}_k\}} \simeq O(k)\times O(n-1-k) \subset O(n-1)$. Using this and the fact that the induced metric on the fibers $Gr_k(T_xX)\simeq Gr(k,n-1)$ is the standard $O(n-1)$-invariant metric we obtain
\begin{equation}
\begin{aligned}
\E_{\Lambda\in Gr(k,n-1)} B_x(\Lambda) &= \frac{1}{|Gr(k,n-1)|}\int\limits_{Gr(k,n-1)} B_x(\Lambda)\, dV_{Gr(k,n-1)} \\
&= \frac{1}{|Gr(k,n-1)|\cdot\left| O(k) \right|\cdot \left| O(n-1-k) \right|}\int\limits_{O(n-1)} \det(E^tg^tDgE)\, dg\\
&=\frac{1}{|O(n-1)|}\int\limits_{O(n-1)} \det(E^tg^tDgE)\, dg
\end{aligned}
\end{equation}
where $dg=dV_{O(n-1)}$ is the invariant Haar measure on $O(n-1)$.

Now for any $k$-subset $I=\{i_1,\dots, i_k\} \subset \{1,\dots,n-1\}$ denote by $E_{I}$ the $(n-1)\times k$ matrix with columns ${\delta}_{i_1},\dots,{\delta}_{i_k}$. $E_I$ can be obtained as a left multiplication of $E$ by the permutation matrix $M_{\sigma_I}$: $E_I = M_{\sigma_I}\cdot E$, where $\sigma_I$ is any permutation that sends $1,\dots,k$ into $i_1,\dots,i_k$ respectively. Using invariance of $dg$ we get
\begin{align}
\int\limits_{O(n-1)}\hspace{-5pt} \det(E_I^t g^t D g E_I)\, dg =\hspace{-5pt} \int\limits_{O(n-1)}\hspace{-5pt} \det(E^t(gM_{\sigma_I})^tD(gM_{\sigma_I})E)\,dg =\hspace{-5pt} \int\limits_{O(n-1)}\hspace{-5pt} \det(E^tg^tDgE)\,dg
\end{align}
Consequently we can express $\E_{\Lambda\in Gr(k,n-1)} B_x(\Lambda)$ as a sum over all $k$-subsets $I \subset \{1,\dots,n-1\}$ divided by ${n-1 \choose k}$:
\begin{align}
\E_{\Lambda\in Gr(k,n-1)} B_x(\Lambda) = {n-1 \choose k}^{-1} \frac{1}{|O(n-1)|} \int\limits_{O(n-1)} \sum\limits_{\substack{I \subset \{1,\dots,n-1\},\\ |I|=k}} \det(E_I^t g^t D g E_I)\, dg
\end{align}
The integrand here is the sum of all principal minors of $g^t D g$ of order $k$ and thus does not depend on $g$ and is equal to the $k$-th elementary symmetric polynomial $s_k(d_1,\dots,d_{n-1})$ of $d_1\,\dots,d_{n-1}$. Combining this with \eqref{eq:step1} and \eqref{inner integral} we end the proof.
\end{proof}
In particular we can derive the following corollary.
\begin{cor}\label{cor:convex}If $X\subset \R\textrm{\emph{P}}^n$ is a convex hypersurface, then
\begin{align}\label{eq:convex} \frac{|\Omega_k(X)|}{|\textrm{\emph{Sch}}(k,n)|}&=\frac{\Gamma\left(\frac{k+1}{2}\right)\Gamma\left(\frac{n-k}{2}\right)}{\pi^{\frac{n+1}{2}}}\int_{X}\sigma_k(x) dV_X.\end{align}
\end{cor}
\begin{proof}We first observe that
\be \frac{|Gr(k,n-1)|}{|\G(k,n)|}=\frac{1}{\pi^{n/2}}\frac{\Gamma\left(\frac{n}{2}\right) \Gamma\left(\frac{n+1}{2}\right)}{\Gamma\left(\frac{k+1}{2}\right)\Gamma\left(\frac{n-k}{2}\right)}\ee
and, recalling \cite[Theorem 4.2]{PSC},
\be \frac{|\textrm{Sch}(k,n)|}{|\G(k,n)|}=\frac{|\Sigma(k+1, n+1)|}{|Gr(k+1, n+1)|}
   =\frac{\Gamma\left(\frac{k+2}{2}\right)}{\Gamma\left(\frac{k+1}{2}\right)} \cdot
      \frac{\Gamma\left(\frac{n-k+1}{2}\right)}{\Gamma\left(\frac{n-k}{2}\right)} .\ee
      Substituting into \eqref{eq:volume convex} we obtain
      \begin{align}
     \frac{|\Omega_k(X)|}{|\textrm{Sch}(k,n)|} &=\frac{|Gr(k,n-1)|}{|\G(k,n)|}\cdot\frac{|\G(k,n)|}{|\textrm{Sch}(k,n)|}\cdot \frac{1}{{n-1 \choose k}}\int_{X}\sigma_k(x) dV_X\\
     &=\frac{1}{\pi^{n/2}}\frac{\Gamma\left(\frac{n}{2}\right) \Gamma\left(\frac{n+1}{2}\right)}{\Gamma\left(\frac{k+1}{2}\right)\Gamma\left(\frac{n-k}{2}\right)}\cdot \frac{\Gamma\left(\frac{k+1}{2}\right) \Gamma\left(\frac{n-k}{2}\right)}{\Gamma\left(\frac{k+2}{2}\right)\Gamma\left(\frac{n-k+1}{2}\right)}\cdot \frac{1}{{n-1 \choose k}}\int_{X}\sigma_k(x) dV_X\\
     &=\frac{1}{\pi^{n/2}}\frac{\Gamma\left(\frac{n}{2}\right) \Gamma\left(\frac{n+1}{2}\right)}{\Gamma\left(\frac{k+2}{2}\right)\Gamma\left(\frac{n-k+1}{2}\right)}\cdot\frac{1}{{n-1 \choose k}}\int_{X}\sigma_k(x) dV_X\\
     &=\frac{\Gamma\left(\frac{k+1}{2}\right)\Gamma\left(\frac{n-k}{2}\right)}{\pi^{\frac{n+1}{2}}}\int_{X}\sigma_k(x) dV_X.
      \end{align}\end{proof}

\subsection{Intrinsic volumes}
Recall that the intrinsic volumes $V_0(C),\dots, V_{n-1}(C)$ of a convex set $C\subset \RP^n$ are characterized by Steiner's formula, which gives the exact expansion (for small $\epsilon>0$) of the volume of the $\epsilon$-neighbourhood of $C$:
 \be\label{eq:steiner} |\mathcal{U}_{\,\RP^n}(C,\epsilon)|=|C|+\sum_{k=0}^{n-1}f_{k}(\epsilon)|S^k||S^{n-k-1}|V_{k}(C)\ee
 (the functions $f_{k}$ are defined in \eqref{eq:universal}). The formula \eqref{eq:steiner} is obtained from the spherical Steiner's formula \cite[(9)]{GaoSchneider} as follows. For a convex set $C\subset\RP^n$ denote by $\tilde C\subset S^n$ any of the two components of $p^{-1}(C)$, where $p: S^n \rightarrow \RP^n$ is the double covering. Under $p$ an open hemisphere in $S^n$ maps isometrically onto $\RP^n$ minus a hyperplane. Therefore, for small $\varepsilon>0$ we have $|\mathcal{U}_{\,\RP^n}(\tilde C,\epsilon)|=|\mathcal{U}_{\,\RP^n}( C,\epsilon)|$ and $V_j(\tilde C) = V_j(C), j=0,\dots,n-1$.
As a consequence we obtain.
\begin{cor}[The manifold of $k$-tangents and intrinsic volumes]Let $C\subset \R\emph{\textrm{P}}^n$ be a strictly convex set with the smooth boundary $\partial C$. Then
\be 4 \cdot V_{n-k-1}(C)=\frac{|\Omega_{k}(\partial C)|}{|\textrm{\emph{Sch}}(k, n)|},\quad k=0, \ldots, n-1.\ee
\end{cor}
\begin{proof}
From \cite[(10)]{GaoSchneider} and Corollary \ref{cor:convex} it follows that
\begin{align} V_{n-k-1}(C)&=\frac{1}{|S^{k}||S^{n-k-1}|}\int_{\partial C}\sigma_k(x)dV_{\partial C}\\
&=\frac{\Gamma\left(\frac{k+1}{2}\right)\Gamma\left(\frac{n-k}{2}\right)}{4\pi^{\frac{n+1}{2}}}\int_{\partial C}\sigma_k(x) dV_{\partial C}\\
&=\frac{1}{4}\cdot \frac{|\Omega_{k}(\partial C)|}{|\textrm{Sch}(k, n)|}.
\end{align}
%We write Spherical Steiner's formula in two different ways. On one hand \cite{GaoSchneider} we have:
%\be \label{eq:tube1}|C_\epsilon|=|C|+\sum_{k=0}^{n-1}|T_{S^n}(S^k, \epsilon)|\cdot V_{k}(C),\ee
%where $T_{S^n}(S^k, \epsilon)$ denotes the $\epsilon$-neighborhood of an equator $S^k$ in $S^n$.
%On the other hand, in the smooth case \cite{GaoSchneider} we have:
%\be \label{eq:tube2} |C_\epsilon|=|C|+\sum_{k=0}^{n-1}|T_{S^n}(S^k, \epsilon)||S^{k}||S^{n-k-1}|\int_{\partial C}\sigma_{n-k-1}(x)dV_{\partial C}.\ee
%Since the function $f_{n, 0}, \ldots, f_{n, n-1}$ are independent, comparing \eqref{eq:tube1} with \eqref{eq:tube2} we obtain
%\begin{align} V_{n-k-1}(C)&=\frac{1}{|S^{k}||S^{n-k-1}|}\int_{\partial C}\sigma_k(x)dV_{\partial C}\\
%&=\frac{\Gamma\left(\frac{k+1}{2}\right)\Gamma\left(\frac{n-k}{2}\right)}{4\pi^{\frac{n+1}{2}}}\int_{X}\sigma_k(x) dV_{\partial C}\\
%&\stackrel{\eqref{eq:convex}}{=}\frac{1}{4}\cdot \frac{|\Omega_{k}(\partial C)|}{|\textrm{Sch}(k, n)|}.
%\end{align}
\end{proof}
This together with \cite[(15)]{GaoSchneider} implies the following interesting corollary.
\begin{cor}\label{cor:sum}Let $C\subset \R\emph{\textrm{P}}^n$ be a strictly convex set with the smooth boundary $\partial C$ and let $C^\circ$ be the polar set of $\tilde C\subset S^n$. Then
\be\frac{4|C|}{|S^n|}+\frac{4|C^\circ|}{|S^n|}+ \sum_{k=0}^{n-1} \frac{|\Omega_k(\partial C)|}{|\textrm{\emph{Sch}}(k,n)|}=4.\ee
In particular, for every $k=0, \ldots, n-1$ we have
\be\label{eq:upperbound}\frac{|\Omega_k(\partial C)|}{|\textrm{\emph{Sch}}(k,n)|}\leq 4.\ee
\end{cor}

\section{Hypersurfaces in random position}
%\subsection{Probabilistic enumerative geometry}

\begin{thm}\label{thm:main}The average number of $k$-planes in $\R{\emph{\textrm{P}}}^n$ simultaneously tangent to convex hypersurfaces  $X_1, \ldots, X_{d_{k,n}}\subset \R\txt{P}^n$ in random position equals
 \be\label{eq:ring} \tau_k(X_1, \ldots, X_{d_{k,n}})={\delta}_{k,n} \cdot \prod_{i=1}^{d_{k,n}}\frac{|\Omega_k(X_i)|}{|\emph{\textrm{Sch}}(k,n)|}.\ee
\end{thm}
\begin{proof}We use the generalized kinematic formula for coisotropic hypersurfaces of $\G(k,n)$ proved in \cite{PSC} (Theorem \ref{Poincare theorem} above). 

In order to apply Theorem \ref{Poincare theorem} to the case $\mathcal{H}_i=\Omega_{k}(X_i), i=1, \ldots, d_{k,n}$,  we need to prove that each $\Omega_{k}(X_i)$ is a  coisotropic hypersurface of $\G(k,n)$. Given $(x,\Lambda) \in Gr_k(X_i)$ as in the proof of Proposition \ref{propo 1} let's consider an orthonormal basis $v_1,\dots,v_{n-1}$ of $T_x X_i$ such that $\Lambda= \textrm{span}\{v_1,\dots,v_k\}$ and a unit normal $\nu \in T_x\RP^n$ to $T_xX_i$. For a curve $x_{\nu}(t)\subset \RP^n$ through $x$ in the direction $\nu$ we consider the parallel transports $v_1^{\nu}(t),\dots,v_k^{\nu}(t) \in T_{x_{\nu}(t)}\RP^n$ of $v_1,\dots,v_k$ along $x_{\nu}(t)$. We claim that the tangent vector to the curve 
$\gamma(t) = x_{\nu}(t)\wedge v_1^{\nu}(t)\wedge \cdots\wedge v_k^{\nu}(t) \in \G(k,n)$ 
is normal to $T_{x\wedge v_1\wedge\cdots\wedge v_k}\Omega_{k}(X_i)$. 
Indeed,
$$\dot{\gamma}(0) = \nu\wedge v_1\wedge\cdots \wedge v_k + \sum\limits_{j=1}^k x\wedge v_1\wedge \cdots\wedge \dot{v}_j^{\nu}(0) \wedge \cdots \wedge v_k = \nu\wedge v_1\wedge \cdots\wedge v_k$$
since $\dot{v}_j^{\nu}(0)= \nabla^{\RP^n}_{\nu} v_j^{\nu} + a_j x = 0 + a_j x$ is proportional to $x$. Now it is elementary to verify that $\dot{\gamma}(0)$ is orthogonal to the tangent space $T_{x\wedge v_1\wedge\cdots\wedge v_k}\Omega_{k}(X_i)$ described in \eqref{propo 1}. Seen as an operator $\dot{\gamma}(0)$ sends $x$ to $\nu$ and all vectors in $\Lambda$ to $0$. Hence $\Omega_{k}(X_i)$ is coisotropic.

Applying now Theorem \ref{Poincare theorem}  we deduce
\be\label{eq:interm} \tau_k(X_1, \ldots, X_{d_{k,n}})=\alpha(k+1,n-k)\,|\G(k,n)|\,\prod\limits_{i=1}^{d_{k,n}}\frac{|\Omega_{k}(X_i)|}{|\G(k,n)|}.\ee
Note that applying Theorem \ref{Poincare theorem} to the special real Schubert variety $\text{\normalfont{Sch}}(k,n)$ we obtain
\begin{align}
{\delta}_{k,n} &= \E \#(g_1\text{\normalfont{Sch}}(k,n)\cap\cdots\cap g_{d_{k,n}}\text{\normalfont{Sch}}(k,n))\\
\label{expected degree}&= \alpha(k+1,n-k)\,|\G(k,n)|\,\left(\frac{|\text{\normalfont{Sch}}(k,n)|}{|\G(k,n)|}\right)^{d_{k,n}}.
\end{align}
This gives an expression for $\alpha(k+1, n-k)$, which substituted into \eqref{eq:interm} gives \eqref{eq:ring}.
\end{proof}
As a consequence we derive the following corollary, which gives a universal upper bound to our random enumerative problem.
\begin{cor}\label{cor:upper}If $X_1, \ldots, X_{d_{k,n}}\subset \R\txt{P}^n$ are convex hypersurfaces, then
\be \tau_k(X_1, \ldots, X_{d_{k,n}})\leq {\delta}_{k,n}\cdot 4^{d_{k,n}}.\ee
\end{cor}
\begin{proof}This follows immediately from \eqref{eq:ring} and \eqref{eq:upperbound}.
\end{proof}
\section{The semialgebraic case}\label{sec:nonconvex}
In this section we discuss a generalization of the previous results to the case of semialgebraic hypersurfaces satisfying some nondegeneracy conditions.

Let $X$ be a smooth closed semialgebraic hypersurface in $\RP^n$. As in Section \ref{sec:convex} we define the Grassmannian bundle of $k$-planes over $X$:
\begin{equation*}
\begin{aligned}
  p: Gr_k(X) &\rightarrow X\\
  (x,\Lambda) &\mapsto x\\
Gr_k(X):&=\{(x,\Lambda):  x\in X, \Lambda\in Gr_k(T_xX)\simeq Gr(k,n-1)\}
\end{aligned}
\end{equation*}
The variety $\Omega_k(X)$ of $k$-tangents to $X$ coincides with the image $\textrm{im}(\psi)$ of the \emph{$k$th Gauss map}:
\begin{align*}
  \psi: Gr_k(X) &\rightarrow \G(k,n)\\
  (x,\Lambda) &\mapsto \text{\normalfont{P}}(\text{\normalfont{Span}}\{x,\Lambda\})
\end{align*}
but now, unlike to the case of a convex hypersurface, $\Omega_k(X)$ is in general singular.

It is convenient to identify the smooth manifold $Gr_k(X)$ with its image in $X\times \G(k,n)$ under the map $\text{\normalfont{id}}\times \psi$:
\begin{align*}
Gr_k(X)\simeq (\text{\normalfont{id}}\times \psi)(Gr_k(X)) &=\{(x,\Lambda)\in X\times \G(k,n): T_x\Lambda \subset T_xX\}\\
 \text{\normalfont{id}}\times \psi: Gr_k(X) &\rightarrow X\times \G(k,n)\\
  (x,\Lambda) &\mapsto (x,\text{\normalfont{P}}(\text{\normalfont{Span}}\{x,\Lambda\}))
\end{align*}
Note that $Gr_k(X)$ is a smooth semialgebraic subvariety of $X\times \G(k,n)$ and the variety of tangents $\Omega_k(X)$ is obtained by projecting it onto the second factor. 

For a point $x\in X$ let's denote by $B$ the second fundamental form of $X$ defined locally near $x$ using any of the two local coorientations of $X$. For $(x,\Lambda)\in Gr_k(X)$ and an orthonormal basis $v_1,\dots,v_k$ of $\Lambda$ denote by $B_x(\Lambda)=\det(B(v_i,v_j))$ the determinant of the $k\times k$ matrix $\{B(v_i,v_j)\}$. Notice that $|B_x(\Lambda)|$ does not depend on the choice of $v_1,\dots,v_k$ and the local coorientation of $X$ near $x$. 
\begin{defi}\label{nd}
We say that $X\subset \R\txt{P}^n$ is $k$-\emph{non-degenerate} if
\begin{enumerate}
\item the semialgebraic set
\begin{align*}
 D:=\{\Lambda \in \G(k,n): \#(\psi^{-1}(\Lambda))>1\}\subset \Omega_k(X)
\end{align*}
of $k$-flats that are tangent to $X$ at more than one point has codimension at least one in $\Omega_k(X)$ and
\item the semialgebraic set
\begin{align*}
  S: = \{ (x,\Lambda) \in Gr_k(X): B|_{T_x\Lambda}\ \text{is degenerate} \} 
\end{align*}
has codimension at least one in the semialgebraic variety $Gr_k(X)$.
\end{enumerate}
\end{defi}
\begin{remark}\label{nd_remark}
Note that the sets $D$ and $S$ are closed in $\Omega_k(X)$ and $Gr_k(X)$ respectively and, by the same reasoning as in the proof of Proposition \ref{propo 1} (up to \eqref{criteria}), the set $S$ consists of such $(x,\Lambda)\in Gr_k(X)$ where $\pi_2: Gr_k(x) \rightarrow \G(k,n)$ is not an immersion.
\end{remark}
A convex semialgebraic hypersurface is $k$-non-degenerate for any $k=0,\dots,n-1$ since in this case the sets $D, S$ from Definition \ref{nd} are empty. The following lemma shows that a generic algebraic surface in $\R\txt{P}^3$ of sufficiently high degree is $1$-non-degenerate.
\begin{lemma}\label{lemma:nondeg}
Let $X^\C\subset \C\txt{P}^3$ be an irreducible smooth surface of degree $d\geq 4$ which does not contain any lines and such that $X=\R X^\C\subset \R\txt{P}^3$ is of dimension $2$. Then $X$ is $1$-non-degenerate.
\end{lemma}
\begin{proof}
Theorem $4.1$ in \cite{kohn:16} asserts that under the assumptions of the current lemma the singular locus $\Sigma^{\C}:=\text{\normalfont{Sing}}(\Omega_1(X^\C))$ of the variety $\Omega_1(X^\C)\subset \G^\C(1,3)$ of complex lines tangent to the complex surface $X^\C\subset \CP^3$ is described as follows:
\begin{align*}
\Sigma^\C = D^\C\cup I^\C,
\end{align*}
where $D^\C$ consists of lines that are tangent to $X^\C$ at more than one point and $I^\C$ consists of lines intersecting $X^\C$ at some point with multiplicity at least $3$. 

We now show that the singular locus $\Sigma:=\txt{Sing}(\Omega_1(X)) = \Omega_1(X)\cap\txt{Sing}(\Omega_1^\C(X))$ of $\Omega_1(X)$ is of dimension at most $2$. There are two cases: either $(1)$ there exists $\Lambda\in \Sigma$ which is smooth for both $\Sigma$ and $\Sigma^\C$ or $(2)$ any smooth point $\Lambda\in \Sigma$ of $\Sigma$ is singular for $\Sigma^\C$. In the case $(1)$ we have $\txt{dim}_{\R}(\Sigma)=\txt{dim}_{\R}(T_{\Lambda}\Sigma)= \txt{dim}_{\C}(T_{\Lambda}\Sigma^\C)=\txt{dim}_{\C}(\Sigma^\C)<\txt{dim}_{\C}(\Omega_1(X^\C)) =3$ and therefore $\txt{dim}_\R(\Sigma)\leq 2$. In the case $(2)$ we have $\txt{dim}_\R(\Sigma) = \txt{dim}_\R(T_{\Lambda}\Sigma) \leq \txt{dim}_\C(\txt{Sing}(\Sigma^\C))< \txt{dim}_\C(\Sigma^\C)< \txt{dim}_{\C}(\Omega_1(X^\C)) =3$ and hence $\txt{dim}_\R(\Sigma)\leq 1$.

For the complex surface $X^\C\subset \RP^3$ let $Gr_1(X^\C)= \{(x,\Lambda)\in X^\C\times \G^\C(1,3): T_{x}\Lambda\subset T_xX^\C\}$ be the Grassmannian bundle of complex lines over $X^\C$.
In the proof of \cite[Thm. 4.1]{kohn:16} it is shown that a line $\Lambda\in I^\C$ intersects $X^\C\subset \CP^3$ at a point $x\in X^\C$ with multiplicity at least $3$ if and only if the differential $(\pi_2)_*: T_{(x,\Lambda)}Gr_1(X^\C) \rightarrow T_{\Lambda}\G^\C(1,3)$ is not injective. By \eqref{criteria} for $(x, \Lambda)\in S$ the differential $(\pi_2)_*: T_{(x,\Lambda)}Gr_1(X) \rightarrow T_{\Lambda}\G(1,3)$ (and hence also $(\pi_2)_*: T_{(x,\Lambda)}Gr_1(X^\C) \rightarrow T_{\Lambda}\G^\C(1,3)$) is not injective. In particular $\pi_2(S)\subset \Omega_1(X)\cap I^\C.$ Now, if $X^\C$ does not contain any lines, the fibers of the projection $\pi_2: Gr_1(X)\rightarrow \G(1,3)$ are finite and hence $\dim(\pi_2(S))=\dim S$. On the other hand, since $\pi_2(S)\subset \Sigma$, the above arguments show that $\dim (\pi_2(S))\leq \txt{dim}(\Sigma)\leq 2$ and consequently $\txt{dim}(S)\leq 2<3=\txt{dim}(Gr_1(X))$. Moreover, this together with \eqref{criteria} imply that there exists a point in $Gr_1(X)$ at which $\pi_2: Gr_1(X)\rightarrow \G(1,3)$ is an immersion and hence $\Omega_1(X)=\pi_2(Gr_1(X))$ is of dimension $3$. 

Observe finally that $D\subset \Omega_1(X)\cap D^\C\subset \Sigma$ and the above arguments imply that $\txt{dim}(D)\leq \txt{dim}(\Sigma)\leq 2<3=\txt{dim}(\Omega_1(X))$.
This finishes the proof.
\end{proof}
\begin{remark}
The above lemma implies that a generic algebraic surface $X\subset \RP^3$ of high enough degree is $1$-non-degenerate.
\end{remark}
In the following proposition we provide a formula for the volume of $\Omega_k(X)$.
\begin{prop}
Let $X$ be a $k$-non-degenerate semialgebraic hypersurface in $\R\text{\normalfont{P}}^n$. Then
\begin{align}\label{volume for nd hypersurfaces}
 |\Omega_k(X)| = |Gr(k,n-1)| \int\limits_{X} \E_{\Lambda\in Gr(k,n-1)} |B_x(\Lambda)|\, dV_X
\end{align}
\end{prop}
\begin{proof}
The complement
\begin{align*}
R := Gr_k(X)\setminus S = \{(x,\Lambda) \in Gr_k(X): B|_{T_x\Lambda}\ \text{\normalfont{is non-degenerate}}\}
\end{align*}
of $S$ is an open dense semialgebraic subset of $Gr_k(X)$. Let's pull back the metric from $\G(k,n)$ to $R$ through the immersion $\pi_2|_R$. Then repeating the proof of Proposition \ref{propo 1} up to the point \eqref{eq:step1} we get
\begin{align}\label{integral over X_R}
  \int\limits_{R} dV_{R} = \int\limits_{X_R} \,\int\limits_{\Lambda \in \pi_1^{-1}(x)\cap R} |B_x(\Lambda)|\, dV_{\pi_1^{-1}(x)\cap R}\,dV_{X_R}
\end{align}
where $X_R:=\pi_1(R)\subset X$ is the projection of $R\subset X\times \G(k,n)$ onto the first factor and the fiber $\pi_1^{-1}(x)=Gr_k(T_xX) \simeq Gr(k,n-1) \subset Gr_k(X)$ is endowed with the uniform distribution.
Note that since $B_x(\Lambda) = 0$ precisely for $\Lambda \in \pi^{-1}_1(x)\setminus R$ we can extend the integration over the whole fiber $\pi^{-1}_1(x)$ in \eqref{integral over X_R}. Moreover, since $X_R=\pi_1(R)$ is open and dense in $X$ (being the image of an open and dense set under the projection $\pi_1$) and since the function
\begin{align*}
  x \mapsto \int\limits_{\Lambda\in\pi_1^{-1}(x)} |B_x(\Lambda)|\, dV_{\pi_1^{-1}(x)}
\end{align*}
is continuous 
\eqref{integral over X_R} becomes
\begin{align}
  \int\limits_{R} dV_{R} = \int\limits_{X} \int\limits_{\Lambda\in\pi_1^{-1}(x)} |B_x(\Lambda)|\, dV_{\pi_1^{-1}(x)}\, dV_X=|Gr(k,n-1)| \int\limits_X \E_{\Lambda\in Gr(k,n-1)} |B_x(\Lambda)|\,dV_X
\end{align}
It remains to prove that $|\Omega_k(X)| = \int\limits_{R} dV_R$. For this let's consider the set 
\begin{align*}
\tilde{D}:=\pi_2^{-1}(D) = \{(x,\Lambda)\in Gr_k(X): \#(\pi_2^{-1}(\Lambda)) > 1\}
\end{align*} 
Note that $\tilde{D}$ is a closed semialgebraic subset of $Gr_k(X)$ and from Definition \ref{nd} it follows that $\tilde{D}\subset Gr_k(X)$ is of codimension at least one. As a consequence, the semialgebraic set $R\setminus \tilde D$ is open and dense in $Gr_k(X)$ (and hence also in $R$) and therefore its projection $\pi_2(R\setminus \tilde{D})$ is open and dense in $\Omega_k(X)$. In particular,
\begin{align*}
  |\Omega_k(X)| = |\pi_2(R\setminus\tilde{D})|=\int\limits_{R\setminus \tilde{D}} dV_{R\setminus \tilde D} = \int\limits_{R} dV_R 
\end{align*}
\end{proof}
\begin{remark}
Using, for example, the Cauchy-Binet theorem it is easy to derive the inequality
 \begin{align*}
    |\Omega_k(X)| \leq \frac{|Gr(k,n-1)|}{{n-1 \choose k}}\int_X s_k(|d_1(x)|,\dots,|d_{n-1}(x)|)\, dV_X
  \end{align*}
where $s_k(|d_1(x)|,\dots,|d_{n-1}(x)|)$ is the $k$th elementary symmetric poynomial of the absolute principal curvatures at $x\in X$. Unfortunately, we do not have a clear geometric interpretation of the right-hand side of the above inequality.

\end{remark}
In the case of lines tangent to a surface in $\R\txt{P}^3$ we can refine the formula \eqref{volume for nd hypersurfaces} as follows.
\begin{cor} If $X \subset \R\txt{P}^3$ is a smooth $1$-non-degenerate surface then
  \begin{align*}
    |\Omega_1(X)| = \int\limits_X h(d_1(x),d_2(x))\,dV_X
  \end{align*}
where 
\begin{align*}
  h(d_1,d_2) = 
\left\{
  \begin{aligned}
    &\frac{\pi}{2}\,|d_1+d_2|,\ \textrm{if}\ d_1d_2\geq 0\\
    & 2\sqrt{-d_1d_2}+2|d_1+d_2|\cdot\left|\arctan\sqrt{-\frac{d_1}{d_2}}-\frac{\pi}{4}\right|,\ \textrm{if}\ d_1d_2<0
  \end{aligned}
\right.
\end{align*}
and $d_1(x),d_2(x)$ are the principal curvatures of $X$ at the point $x$.
\end{cor}
\begin{proof}
The formula \eqref{volume for nd hypersurfaces} reads
\begin{align*}
  |\Omega_1(X)| = \pi \int\limits_X \E_{\Lambda \in Gr(1,2)} |B_x(\Lambda)|\, dV_X
\end{align*}
In coordinates in which the second fundamental form $B_x$ of $X\subset \RP^3$ at the point $x\in X$ is diagonal with values $d_1,d_2$ we have
\begin{align*}
\pi\,  \E_{\Lambda\in Gr(1,2)} |B_x(\Lambda)| =\pi\, \E_{v\in S^1} |B_x(v,v)| = \int_{-\pi/2}^{\pi/2} |d_1\cos^2 \varphi + d_2\sin^2 \varphi| d\varphi 
\end{align*}
The last integral can be evaluated by elementary integration methods giving $\frac{\pi}{2}|d_1+d_2|$ in case $d_1d_2\geq 0$ and
\begin{align*}
  2\sqrt{-d_1d_2}+2|d_1+d_2|\cdot\left|\arctan\sqrt{-\frac{d_1}{d_2}}-\frac{\pi}{4}\right|
\end{align*}
in case $d_1d_2<0$. 
\end{proof}

Finally we prove an analog of Theorem \ref{thm:main} for $k$-non-degenerate semialgebraic hypersurfaces.
\begin{thm}\label{thm:hyper}
The average number of $k$-flats in $\R{\emph{\textrm{P}}}^n$ simultaneously tangent to $k$-non-degenerate semialgebraic hypersurfaces  $X_1, \ldots, X_{d_{k,n}}$ in random position equals
 \be \tau_k(X_1, \ldots, X_{d_{k,n}})={\delta}_{k,n} \cdot \prod_{i=1}^{d_{k,n}}\frac{|\Omega_k(X_i)|}{|\emph{\textrm{Sch}}(k,n)|}.\ee
\end{thm}
\begin{proof}
Exactly in the same way as in the proof of \ref{thm:main} one can show that the smooth locus $\Omega_k(X_i)_{sm}$ of $\Omega_k(X_i)$ is a coisotropic hypersurface of $\G(k,n)$. Since $\Omega_k(X_i)\setminus \Omega_k(X_i)_{sm}$ has codimension $\geq 2$ in $\G(k,n)$ by standard transversality arguments we have that  
\begin{align*}
g_1\Omega_k(X_1) \cap \cdots \cap g_{d_{k,n}}\Omega_k(X_{d_{k,n}}) = g_1\Omega_k(X_1)_{sm} \cap \cdots \cap g_{d_{k,n}} \Omega_k(X_{d_{k,n}})_{sm}
\end{align*}
for a generic choice of $g_1,\dots,g_{d_{k,n}} \in O(n+1)$.

The claim follows by applying the integral geometry formula (Theorem \ref{Poincare theorem}) to the semialgebraic sets $\Omega_k(X_1)_{sm},\dots,\Omega_k(X_{d_{k,n}})_{sm}$ as in the proof of Theorem \ref{thm:main}.
\end{proof}
\begin{remark}(Random invariant hypersurfaces)
The previous Theorem can be used for computing the expectation of the number of $k$-flats tangent to random Kostlan hypersurfaces of degree $m_1, \ldots, m_{d_{k,n}}$ in $\RP^n$ -- notice that here the randomness comes directly from the hypersurfaces! Let us discuss the case $n=3$, $k=1$. 

Let $f_1, \ldots, f_4\in \R[x_1,\dots,x_4]$ be random, independent, $O(4)$-invariant polynomials of degree $m_1, \ldots, m_4\geq 4$ and denote by $X(f_i)=\{f_i=0\}\subset \RP^3$, $i=1, \ldots, 4$ the corresponding projective hypersurfaces. We are interested in computing
\be (*)=\E_{f_1, \ldots, f_4}\#\Omega_1(X(f_1))\cap\cdots\cap\Omega_1(X(f_4)).\ee
We use the fact that the polynomials are invariant for writing:
\begin{align}(*)&=\E_{g_1, \ldots, g_4}\E_{f_1, \ldots, f_4}\#\Omega_1(g_1X(f_1))\cap\cdots\cap\Omega_1(g_4X(f_4))\\
&=\E_{f_1, \ldots, f_4}\E_{g_1, \ldots, g_4}\#\Omega_1(g_1X(f_1))\cap\cdots\cap\Omega_1(g_4X(f_4)).
\end{align}
For $i=1, \ldots, 4$ with probability one $X(f_i)$ is irreducible and there are no lines on it; hence by Lemma \ref{lemma:nondeg} with probability one each $X(f_i)$ is $1$-non-degenerate. Applying Theorem \ref{thm:hyper} we conclude that
\be \E_{f_1, \ldots, f_4}\#\Omega_1(X(f_1))\cap\cdots\cap\Omega_1(X(f_4))= \delta_{1,3}\cdot \prod_{i=1}^{4}\frac{\E_{f_i}|\Omega_1(X(f_i))|}{|\textrm{Sch}(1,3)|}.\ee
\end{remark}
\section{Convex bodies with many common tangents}\label{sec:counter}
The purpose of this section is to show that for every $m>0$ there exist convex surfaces $X_1, \ldots, X_4\subset\RP^3$ in general position such that the intersection $\Omega_1(X_1)\cap\cdots\cap \Omega_1(X_4)\subset \G(1,3)$ is transverse and consists of at least $m$ points. We owe the main idea for this to T. Theobald.
\subsection{A coordinate system}\label{sec:coordinate}
\begin{figure}
	\centering
	\includegraphics[scale= 0.7]{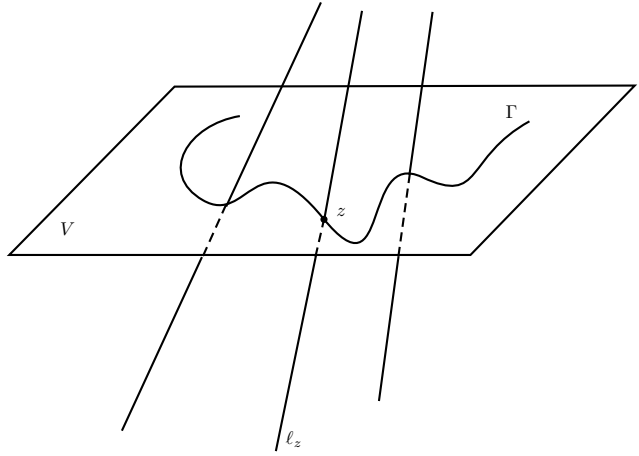}
	\caption{The construction of the coordinate system introduced in Section \ref{sec:coordinate}.} \label{fig:curve1}
\end{figure}
Let $X_1, X_2, X_3\subset\RP^3$ be smooth convex semialgebraic surfaces such that the intersection $Z=\Omega_1(X_1)\cap\Omega_1(X_2)\cap \Omega_1(X_3)$ is transverse (hence $Z$ is a smooth curve in $\G(1,3)$). Let 
\begin{align*}
P = \{ (\Lambda, [v]):\, \Lambda \in Z, [v] \in \Lambda\simeq \RP^1\}
\end{align*} be the projectivized tautological bundle over $Z$ and consider the tautological map
\begin{align*} 
\eta: P &\rightarrow \RP^3\\
  (\Lambda,[v]) &\mapsto [v]
\end{align*}
We determine points where $\eta$ is an immersion.
\begin{lemma}
$\eta_*: T_{(\Lambda, [v])} P \rightarrow T_{[v]}\R\txt{P}^3 \simeq v^{\perp}$ is injective if and only if $v$ is not annihilated by the generator of $T_{\Lambda}Z \subset \text{\normalfont{Hom}}(\Lambda, \Lambda^{\perp})$.  
\end{lemma}
\begin{proof}
  Let $\Lambda(t) = v(t)\wedge u(t)$ be a local parametrization of $Z$ near $\Lambda=\Lambda(0)$, where $\{u(t),v(t)\}$ is an orthonormal basis of $\Lambda(t)$ and $v=v(0), u=u(0)$. The tangent vectors to the curves $\gamma_1(t) = (\Lambda, [\cos t\, v + \sin t\, u])$, $\gamma_2(t) = (\Lambda(t), [v(t)])$ at $t=0$ span the tangent space $T_{(\Lambda,[v])}P$ and $\eta_*(\dot\gamma_1(0)) = [u]$, $\eta_*(\dot \gamma_2(0)) = [\dot v(0)]$. Any generator of the one-dimensional space $T_{\Lambda}Z\subset \text{\normalfont{Hom}}(\Lambda, \Lambda^{\perp})$ sends $v\in \Lambda$ to $\dot v(0) \in \Lambda^{\perp}\subset v^{\perp}$. The assertion follows.
\end{proof}
Let $(\Lambda,[v]) \in P$ be a point where $\eta$ is an immersion (by the above lemma such $(\Lambda,[v])\in P$ exists for any $\Lambda \in Z$) and let $V\simeq \RP^2 \subset \RP^3$ be a plane through $[v]=\eta((\Lambda,[v])) \in \RP^3$ that is transversal to the line $\ell_{[v]} := \eta((\Lambda,\Lambda))$. The map $\eta$ is an embedding locally near $(\Lambda,[v])$. Therefore the image under $\eta$ of a small neighbourhood of $(\Lambda,[v])$ intersects $V$ along a smooth curve which we denote by $\Gamma$. Moreover, the images of the fibers of $P$ define a smooth field of directions $\{\ell_z: z\in \Gamma\}$ on $\Gamma$ (see figure \ref{fig:curve1}) which can be smoothly extended to a field of directions $\{\ell_z: z\in U\}$ on a neighbourhood $U\subset V$ of $\Gamma$.

As a consequence there exists a neighborhood $W\subset \RP^3$ of $[v]$ of the form  
\begin{align*}
W = \coprod\limits_{z\in U} \ell_z \cap W \simeq U\times (-1,1)
\end{align*}
On this neighbourhood we have a smooth map (the projection on the first factor):
\be \pi:W\to U.\ee
This map has the following property:
\begin{lemma}\label{last}
 If $B\subset W$ is a smooth strictly convex subset in $\R\txt{P}^3$ and $z\in U$ is a critical value for $\pi|_{\partial B}$, then $\ell_z$ is tangent to $\partial B$. 
\end{lemma}
\begin{proof}
In fact if $\#\{\ell_z\cap \partial B\}=2$ then the line $\ell_z$ would be trasversal to $\partial B$ and $z$ would be a regular value for $\pi|_{\partial B}$.
\end{proof}

\subsection{The construction}\begin{figure}
	\centering
	\includegraphics[scale= 0.7]{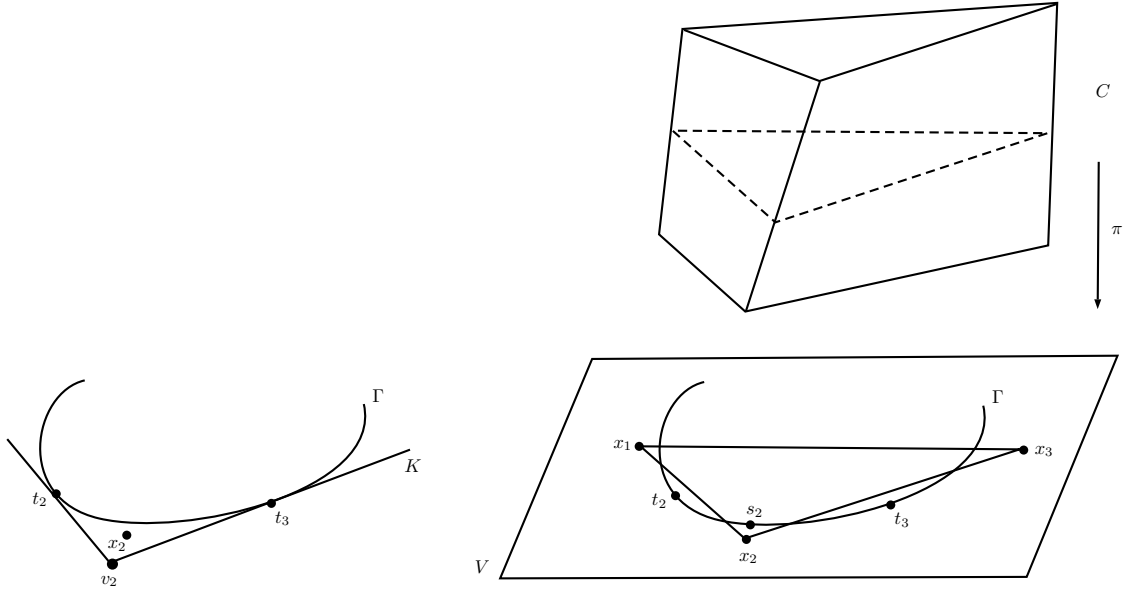}
	\caption{The  convex body $C$.} \label{fig:curve3}
\end{figure}
Using strict convexity of $X_1, X_2, X_3$ it is easy to show that for a generic choice of the plane $V$ a small arc of the curve $\Gamma$ is strictly convex. Let's use the same letter $\Gamma$ to denote such an arc. For a given number $m>0$ pick $n=m+1$ distinct points $t_1, \ldots, t_n$ on $\Gamma$ and consider an $n-$polygonal arc $K$ tangent to $\Gamma$ at the points $t_1, \ldots , t_n$. Call $v_1 \ldots, v_{n-1}$ the ordered vertices of $K$ and  for every (curvilater) triangle $t_iv_it_{i+1}$ pick a point $x_i$ in its interior (see left picture in Figure \ref{fig:curve3}). 

Let now $C\subset W$ be the convex body in $\RP^3$ defined as the convex hull of the segments in $W \simeq U\times (-1,1)$:
\be C=\textrm{conv}(\{x_1\}\times (-\delta, \delta), \ldots, \{x_{n-1}\}\times (-\delta, \delta)),\ee
where $\delta>0$ is chosen small enough such that none of $t_1, \ldots, t_n$ belongs to $\pi(C)$. Note that the polygon $x_1\cdots x_{n-1}$ is a subset of $\pi(C)\subset C$. As a consequence, there exist points $s_1, \ldots, s_{n-1}$ on $\Gamma$, interlacing $t_1, \ldots, t_n$ such that they all belong to $\textrm{im}(\pi|_{\textrm{int}(C)}).$ (See the right picture in Figure \ref{fig:curve3}.)

Let now $C_\epsilon\subset W$ be a smooth, strictly convex semialgebraic approximation of $C$ such that:
\begin{itemize}
	\item [(1)]$s_1, \ldots , s_{n}\in \pi|_{\textrm{int}(C_\epsilon)};$
	\item[(2)] $t_1, \ldots, t_n\notin \pi(C_\epsilon);$
	\item[(3)] the intersection $\Omega_1(C_\epsilon)\cap \Omega_1(X_1)\cap \Omega_1(X_2)\cap \Omega_1(X_3)$ is transverse.
\end{itemize}
The conditions (1) and (2) imply that $\pi(\partial C_\epsilon)\cap \Gamma$ (a semialgebraic subset of $\Gamma$) consists of intervals:
\be \pi(\partial C_\epsilon)\cap \Gamma=[a_1, b_1]\cup \cdots \cup [a_N, b_N]\ee
possibly reduced to points and $N\geq n-1$.
Now each $a_i$ (and $b_i$) is critical for $\pi|_{\partial C_\epsilon}$: otherwise the image of $\pi|_{\partial C_\epsilon}$ near $a_i$ would contain an open set and $a_i$ would not be a boundary point of the intersection $\pi(\partial C_\epsilon)\cap \Gamma$. By Lemma \ref{last} this implies that each line $\ell_{a_i}$ is tangent to $\partial C_\epsilon$ and condition $(3)$ implies that the transverse intersection $\Omega_1(C_\epsilon)\cap \Omega_1(X_1)\cap \Omega_1(X_2)\cap \Omega_1(X_3)$ (which is finite) contains more than $n-1=m$ lines.

\bibliographystyle{plain}
\bibliography{flatsArXiv2}

\end{document}